\documentclass [12pt]  {article} 
\usepackage{amssymb}
\def \NN{{\mathbb{N}}}

\def \FF{{\mathbb{F}}}

\def \tn{{\hbox{$\not|\,$}}}
\def \GGG{{\mathcal{G}}}
\def \UUU{{\mathcal{U}}}

\begin{document}

\begin{center}
{\Large {\bf Several types of solvable groups as automorphism 
groups of compact Riemann surfaces}}\\
\bigskip
\bigskip
{\sc Andreas Schweizer}\\
\bigskip
{\small {\rm Department of Mathematics,\\
Korea Advanced Institute of Science and Technology (KAIST),\\ 
Daejeon 305-701\\
South Korea\\
e-mail: schweizer@kaist.ac.kr}}
\end{center}
\begin{abstract}
\noindent
Let $X$ be a compact Riemann surface of genus $g\geq 2$. Let
$Aut(X)$ be its group of automorphisms and  $G\subseteq Aut(X)$ 
a subgroup. Sharp upper bounds for $|G|$ in terms of $g$ are 
known if $G$ belongs to certain classes of groups, e.g. solvable, 
supersolvable, nilpotent, metabelian, metacyclic, abelian, cyclic.
We refine these results by finding similar bounds for groups of 
odd order that are of these types.
We also add more types of solvable groups to that long list by
establishing the optimal bounds for, among others, groups of 
order $p^m q^n$. 
Moreover, we show that Zomorrodian's bound for $p$-groups $G$ 
with $p\geq 5$, namely $|G|\leq \frac{2p}{p-3}(g-1)$, actually 
holds for any group $G$ for which $p\geq 5$ is the smallest 
prime divisor of $|G|$.
\end{abstract}

\subsection*{1. Introduction}

\noindent
We write $|G|$ for the order of a group $G$ and $G'$ for the
commutator group. Also $C_n$ stands for a cyclic group of order $n$.
\par
In this paper $X$ will always be a compact Riemann surface of genus
$g\geq 2$. Its full group of conformal automorphisms is denoted by
$Aut(X)$.
\par
It is a classical theorem by Hurwitz that then 
$|Aut(X)|\leq 84(g-1).$
See for example [B, Theorem 3.17].
\par
However, not every group whose order is divisible by $84$ can occur
in Hurwitz's Theorem. For example, if $g=7^n +1$, a group of order
$84(g-1)=12\cdot 7^{n+1}$ must be solvable by the Sylow Theorems; 
and for solvable groups there are stronger bounds 
(see Theorem 2.2 (a) below).
\par
More generally, if $G$ is a (not necessarily proper) subgroup 
of $Aut(X)$ and $G$ belongs to a certain type of groups, one 
is interested in having a bound on $|G|$ in terms of $g$.
Ideally one would like to have for any given type of groups 
a simple function $b(g)$ such that $|G|\leq b(g)$ except for 
finitely many (explicitly known) $g$, and that for infinitely 
many values of $g$ there is a $G$ of the desired type with 
$|G|=b(g)$.
\par
Such a function $b(g)$ is of course not guaranteed to exist. 
But for many interesting types of groups it does. From the
rich literature we summarize the results that are relevant 
to this paper.
$$\hbox{\rm TABLE\ 1}$$
$$\def\abstand{$\vrule width 0pt height 18pt $}
\\ \\
\begin{array}{|c|c|c|c|} \hline
\abstand \hbox{\rm  class of $G$} & \hbox{\rm upper bound} & 
\hbox{\rm exceptions} & \hbox{\rm source} \\ 
\hline
\abstand \hbox{\rm general} & 84(g-1) & \hbox{\rm none} & 
\hbox{\rm classical} \\
\abstand \hbox{\rm solvable} & 48(g-1) & \hbox{\rm none} & 
\hbox{\rm [Ch], [G1], [G2]} \\
\abstand \hbox{\rm supersolvable} & 18(g-1) & \hbox{for $g=2$} & 
\hbox{\rm [Z3], [GMl], [Z4]} \\
\abstand \hbox{\rm nilpotent} & 16(g-1) & \hbox{\rm none} & 
\hbox{\rm [Z1]} \\
\abstand \hbox{\rm metabelian} & 16(g-1) & \hbox{\rm for $g=2,3,5$} & 
\hbox{\rm [ChP], [G2], [G3]} \\
\abstand \hbox{\rm metacyclic} & 12(g-1) & \hbox{\rm for $g=2$} & 
\hbox{\rm [Sch2]} \\
\abstand \hbox{\rm $Z$-group} & 10(g-1) & \hbox{\rm for $g=2,3$} & 
\hbox{\rm [Sch2]} \\
\abstand \hbox{\rm $|G|$ is square-free} & 10(g-1) & \hbox{\rm for $g=3$} & 
\hbox{\rm [Sch2]} \\
\abstand \hbox{\rm abelian} & 4g+4 & \hbox{\rm none} & 
\hbox{\rm classical, see also [Ml], [G1]} \\
\abstand \hbox{\rm cyclic} & 4g+2 & \hbox{\rm none} & 
\hbox{\rm classical, see also [Ha], [G1], [N]} \\
\hline
\end{array} $$ 
\\ \\
In Theorem 2.2 in the next section we will report more details
on some of the groups. Here we only mention that the first proof 
for metabelian groups is in [ChP]. But the last theorem of [G2] 
points out that a metabelian group of order $24$ for $g=2$, 
covered by $\Gamma(0;2,4,6)$, had been overlooked. See also [G3] 
for a detailed description of the metabelian $G$ with $|G|=16(g-1)$ 
and some corrections concerning their possible orders.
\par
Bounds and more information if $G$ is a $p$-group were worked out 
in [Z2]; see Theorem 2.3. More recent results in [W2] and [MZ2] 
include the optimal bound $b(g)$ for subgroups $G$ of odd order 
in $Aut(X)$, which we recall in Theorem 2.4.
\\ \\
In this paper we present several new results of a similar 
nature.
\par
In Section 5 we refine the known results by determining 
for each type of group in Table 1 the sharp bound $b(g)$ 
if $G\subseteq Aut(X)$ is an {\it odd order} subgroup of 
that type. It turns out that odd order supersolvable $G$ 
of maximal possible order are subject to much stronger 
restrictions than the other types (see Theorem 5.1). 
\par
The idea of considering groups of odd order can be generalized 
to prescribing the smallest prime divisor $p$ of $|G|$. We do
this already in Section 4, before treating odd order groups,
for reasons having to do with the logical dependence of the 
proofs.
It turns out that for $p\ge 5$ the bounds for the first 8 types 
of groups in Table 1 all collapse to the same bound, namely the 
one for $p$-groups from Theorem 2.3 (c).
\par
For a prime number $p$ denote by $\GGG(p)$ the class of all 
finite groups whose orders are not divisible by any primes 
smaller than $p$. So $\GGG(2)$ denotes all finite groups and 
$\GGG(3)$ all groups of odd order.
Omitting some of the details, the following big picture emerges.
\par
The bound for $p$-groups from Theorem 2.3 is also the bound for 
nilpotent groups inside the class $\GGG(p)$, and any nilpotent
group in $\GGG(p)$ that reaches this bound must be a $p$-group.
Moreover, with the exception of $6$ individual groups ($3$ in
$\GGG(2)$ and $3$ in $\GGG(3)$), this bound also is the sharp 
bound for metabelian groups in $\GGG(p)$. 
If $p\geq 3$ and $q\equiv 1\ mod\ p$ are primes, we can even 
find infinitely many metabelian, supersolvable $(p,q)$-groups 
that attain this bound. 
\par
So in Section 6 we treat yet another type of groups that are known 
to be solvable, namely groups of order $p^m q^n$ where $p<q$ 
are primes. We use the results from the previous sections to get 
optimal bounds on $|G|$ for $G\subseteq Aut(X)$ with $|G|=p^m q^n$, 
in general and if $p$ and/or $q$ is prescribed.
\par
In Section 7 we get similar results for two (probably less important)
types of groups that are lying strictly between the supersolvable and 
the solvable ones, namely groups with nilpotent commutator group, and
groups in which the elements of odd order form a (necessarily normal)
subgroup.
\par
In Section 8 we derive sharp upper bounds on $|G|$ for CLT groups
of odd order. CLT groups are another type of groups between the 
supersolvable and the solvable ones. They are defined in elementary 
terms by the condition that for every divisor $d$ of $|G|$ there 
is a subgroup of order $d$. But CLT groups are not at all 
well-behaved. See Section 8 for details. This makes them a very 
unwieldy object to handle and explains perhaps why we only get
partial results for the size of CLT groups in general.
\\

\subsection*{2. Known results}

\noindent
Let $X$ be a compact Riemann surface of genus $g\geq 2$ and let 
$G$ be a subgroup of $Aut(X)$. The key tool to get information 
on $G$ is that $G\subseteq Aut(X)$ is covered by a Fuchsian
group $\Gamma=\Gamma(h;m_1,\ldots,m_r)$. This means that 
$\Gamma$ is discretely embedded into $Aut(\UUU)$, where $\UUU$
is the complex upper halfplane, and that there is a torsion-free
normal subgroup $K$ of $\Gamma$ with $\Gamma/K\cong G$ such that 
$\UUU/K\cong X$ is the universal covering. Moreover, $h$ is the 
genus of $X/G\cong\UUU/\Gamma$. See Section 3 in Chapter 1 of 
[B] for more background and references. 
\\ \\
{\bf Theorem 2.1.} \it 
If $G$ is covered by a Fuchsian group 
$\Gamma(h;m_1 ,m_2 ,\ldots, m_r)$, then 
$$|G|=\frac{2}{2h-2+\sum_{i=1}^r (1-\frac{1}{m_i})}(g-1).$$

\rm
\noindent
{\bf Proof.} \rm 
[B, Theorem 3.5 and p.15].
\hfill$\Box$
\\ \\
In particular, if $|G|>4(g-1)$, then $h=0$. In practically 
all cases we are interested in, $\Gamma$ even is a triangle 
group, i.e. $h=0$ and $r=3$. 
\par
Since we are interested in solvable groups $G\subseteq Aut(X)$,
the biggest abelian quotient $G/G'$ is very important for us.
It must be a quotient of $\Gamma/\Gamma'$, the structure of 
which can be easily read off (at least if $h=0$) from the 
generators and relations
$$\Gamma(0;m_1,\ldots,m_r)\cong
\langle x_1 , \ldots, x_r\ |\ x_1^{m_1}=x_2^{m_2}
=\ldots=x_r^{m_r}=x_1 x_2\cdots x_r =1\rangle.$$
For frequent use throughout the paper we take from Table 4.1 in
[GMl] all possible orders $|G|\geq 18(g-1)$ for solvable groups 
$G\subseteq Aut(X)$ together with the corresponding triangle
groups $\Gamma$. 
$$ \hbox{\rm TABLE\ 2}$$
$$\def\abstand{$\vrule width 0pt height 15pt $}
\\ \\
\begin{array}{|c|c|c|c|} \hline
 |G| & \Gamma & \Gamma/\Gamma' & \Gamma' \\ 
\hline
\abstand 48(g-1) & \Gamma(0;2,3,8) & C_2 & \Gamma(0;3,3,4) \\
\abstand 40(g-1) & \Gamma(0;2,4,5) & C_2 & \Gamma(0;5,5,2) \\
\abstand 36(g-1) & \Gamma(0;2,3,9) & C_3 & \Gamma(0;2,2,2,3) \\
\abstand 30(g-1) & \Gamma(0;2,3,10) & C_2 & \Gamma(0;3,3,5) \\
\abstand 24(g-1) & \Gamma(0;2,3,12) & C_6 & \Gamma(1;2) \\
\abstand 24(g-1) & \Gamma(0;2,4,6) & C_2 \times C_2 & \Gamma(0;2,2,3,3) \\
\abstand 24(g-1) & \Gamma(0;3,3,4) & C_3 & \Gamma(0;4,4,4) \\
\abstand 21(g-1) & \Gamma(0;2,3,14) & C_2 & \Gamma(0;3,3,7) \\
\abstand 20(g-1) & \Gamma(0;2,3,15) & C_3 & \Gamma(0;2,2,2,5) \\
\abstand 20(g-1) & \Gamma(0;2,5,5) & C_5 & \Gamma(0;2,2,2,2,2) \\
\abstand \frac{96}{5}(g-1) & \Gamma(0;2,3,16) & C_2 & \Gamma(0;3,3,8) \\
\abstand \frac{56}{3}(g-1) & \Gamma(0;2,4,7) & C_2 & \Gamma(0;7,7,2) \\
\abstand 18(g-1) & \Gamma(0;2,3,18) & C_6 & \Gamma(1;3) \\
\hline
\end{array} $$ 
\\ \\
Actually, Table 4.1 in [GMl] contains four more triangle groups
with $|G|\geq 18(g-1)$, namely $\Gamma=\Gamma(0;2,3,p)$ with 
$p=7, 11, 13, 17$. But for them $\Gamma'=\Gamma$ holds, and 
hence they cannot cover a solvable group. $\Gamma(0;2,3,7)$ 
is of course the triangle group that covers Hurwitz groups.
\par
Once one has a Riemann surface $X$ of genus $g\geq 2$ and 
a subgroup $G$ of $Aut(X)$ with desired properties and size,
one can try to construct infinitely many more examples (with
growing genus) from it. This is usually done by an approach
that goes back to Macbeath [Mb] and is used, among others,
in [Ch], [ChP], [G1], [MZ2] and [W2]. For every natural number 
$n$ there exists a compact Riemann surface $Y$ that 
is a totally unramified Galois cover of $X$ with 
$Gal(Y/X)\cong (C_n)^{2g}$. So the genus of $Y$ is 
$n^{2g}(g-1)+1$. The key point is that $G$ lifts to 
$Y$ in the sense that $Aut(Y)$ has a subgroup $\widetilde{G}$
with $Gal(Y/X)\triangleleft\widetilde{G}$ and 
$\widetilde{G}/Gal(Y/X)\cong G$.
\par
By construction, a property like being solvable is obviously 
passed on from $G$ to $\widetilde{G}$. But more special 
properties, for example being supersolvable, might get 
lost or only hold under certain circumstances.
\par
Now we have most of the tools ready to find sharp bounds for 
$|G|$ in terms of $g$ when $G\subseteq Aut(X)$ belongs to 
a certain class of groups. Table 1 in the Introduction gives 
a selective overview over the rich literature. Here we only 
spell out those cases of which we need details later on.
\\ \\
{\bf Theorem 2.2.} \it 
Let $X$ be a compact Riemann surface of genus $g\geq 2$, and 
let $G$ be a subgroup of $Aut(X)$.
\begin{itemize}
\item[(a)] If $G$ is solvable, then $|G|\leq 48(g-1)$. If $G$ 
attains this bound, then necessarily $G=Aut(X)$. There are 
infinitely many $g$ for which there exists a Riemann surface
$X$ of genus $g$ such that $Aut(X)$ is solvable and of order 
$48(g-1)$.
\item[(b)] If $G$ is supersolvable, then for $g\geq 3$ we have
$|G|\leq 18(g-1)$. Moreover, if $g\geq 3$ the necessary and 
sufficient condition for the existence of a Riemann surface 
$X$ of genus $g$ and a supersolvable $G\subseteq Aut(X)$ with 
$|G|=18(g-1)$ is that $g-1$ is divisible by $9$ and has no 
prime divisors that are congruent to $2$ modulo $3$.
\item[(c)] If $G$ is nilpotent, then $|G|\leq 16(g-1)$. Any
nilpotent $G$ that reaches this bound must be a $2$-group.
\end{itemize}
\rm

\noindent
{\bf Proof.} \rm 
(a) From [B, Lemma 3.18] we get the following two facts: 
Groups $G$ with $|G|=84(g-1)$ have $G'=G$, and hence they are 
not solvable. The next possible order is $48(g-1)$. This also
shows that $|G|=48(g-1)$ implies $G=Aut(X)$, because the only
possible bigger order, $84(g-1)$, is not a multiple. Note
however, that $|Aut(X)|=48(g-1)$ does not automatically imply
solvability.
\par
Chetiya [Ch] has constructed for each $n>0$ a Riemann surface 
of genus $2n^6 +1$ whose automorphism group is solvable of 
order $96n^6$. This is the special case $m=4$ of the more 
general result in [Ch, Theorem 3.2]. The principal idea is 
to apply Macbeath's construction to a Riemann surface of 
genus $3$ whose automorphism group has order $96$. 
By the same method, but starting with the unique Riemann surface
of genus $2$ with automorphism group of order $48$, Gromadzki
[G1, Section 5] for every $n>0$ gets a Riemann surface of genus 
$n^4 +1$ with a solvable automorphism group of order $48n^4$.
\par
See also [G2] for analogous results for groups of solvable 
length $\leq 3$.
\par
(b) [Z3], [GMl] and [Z4]. In [Z3] the possibility that $g-1$ 
might have prime divisors that are congruent to $1$ modulo $3$
is erroneously excluded.
\par
(c) [Z1, Theorems 1.8.4 and 2.1.2]
\hfill$\Box$
\\ \\
To complete Theorem 2.2 (c), and for use in later sections, we state
\\ \\
{\bf Theorem 2.3.} [Z2] \it 
Let $X$ be a compact Riemann surface of genus $g\geq 2$, and let $G$
be a subgroup of $Aut(X)$. 
\begin{itemize}
\item[(a)] If $G$ is a $2$-group, then $|G|\leq 16(g-1)$. 
Moreover, for every $n\geq 4$ there exists a group $G$ of 
order $2^n$ that reaches this bound.
\item[(b)] If $G$ is a $3$-group, then $|G|\leq 9(g-1)$. 
Moreover, for every $n\geq 4$ there exists a group $G$ of 
order $3^n$ that reaches this bound.
\item[(c)] If $G$ is a $p$-group with $p\geq 5$, then 
$|G|\leq\frac{2p}{p-3}(g-1)$. 
Moreover, for every $n\geq 1$ there exists a group $G$ of 
order $p^n$ that reaches this bound.
\end{itemize}
\rm
\bigskip
\noindent
{\bf Proof.} \rm
Theorems 1.1.2, 1.2.1, 1.3.1, and 2.0.1 in [Z2].
\hfill$\Box$
\\ \\
The four biggest possible orders if $G\subseteq Aut(X)$ has 
odd order are given in [MZ2, Proposition 1] together with the 
corresponding triangle groups. For further use throughout the 
paper we list them here. 
$$\hbox{\rm TABLE\ 3}$$
$$\def\abstand{$\vrule width 0pt height 15pt $}
\\ \\
\begin{array}{|c|c|c|c|} \hline
 |G| & \Gamma & \Gamma/\Gamma' & \Gamma' \\ 
\hline
\abstand 15(g-1)           & \Gamma(0;3,3,5)  & C_3 & \Gamma(0;5,5,5) \\
\abstand \frac{21}{2}(g-1) & \Gamma(0;3,3,7)  & C_3 & \Gamma(0;7,7,7) \\
\abstand 9(g-1) & \Gamma(0;3,3,9) & C_3 \times C_3 & \Gamma(1;3,3,3) \\
\abstand \frac{33}{4}(g-1) & \Gamma(0;3,3,11) & C_3 & \Gamma(0;11,11,11) \\
\hline
\end{array} $$ 
\\ \\
See [MZ2] for more information. For example on page 328 of 
[MZ2] Macbeath's method is used to show that there are $G$ 
with $|G|=15(g-1)$ for every $g=5n^{12}+1$.
\par
A more precise upper bound for $G$ of odd order is given in 
[W2, Main Theorem]. It also takes into account the highest 
power of $2$ dividing $g-1$. If $g-1$ is odd, it coincides 
with $15(g-1)$, and otherwise it is smaller.
\par 
Correspondingly, [W2, Theorem 4.2] constructs infinite series 
of such groups. Specializing it to $l=0$ we have $N_l =5$ and 
$h=6$ and get the same infinite series with $|G|=15(g-1)$.
\par
We summarize the most important facts.
\\ \\
{\bf Theorem 2.4.} [W2], [MZ2] \it
Let $X$ be a compact Riemann surface of genus $g\geq 2$, and 
let $G$ be a subgroup of $Aut(X)$. If $|G|$ is odd, then
$|G|\leq 15(g-1)$. 
\par
Conversely, for every odd integer $n$ there exists a compact
Riemann surface $X$ of genus $g=5n^{12}+1$ such that $Aut(X)$
contains a subgroup $G$ of order $15(g-1)$.
\rm
\\ \\
Of course, one should keep in mind throughout this paper that
groups of odd order are always solvable. This is proved in the
monumental article [FT]. So Theorem 2.4 should be compared to
Theorem 2.2 (a) and not to the Hurwitz bound. 
\\

\subsection*{3. Some group theoretic tools}

\noindent
What makes supersolvable groups much more convenient to handle 
than the merely solvable ones is the following fact, that will 
be used repeatedly in the paper.
\\ \\
{\bf Theorem 3.1.} (Zappa's Theorem) \it
Let $|G|=\prod_{i=1}^{s}p_i$ with prime numbers 
$p_1\leq p_2\leq \ldots\leq p_s$.
If $G$ is supersolvable, then there exist normal subgroups
$G_i$ of $G$ with
$$G=G_0 \triangleright G_1 \triangleright \ldots \triangleright 
G_{s-1} \triangleright G_s=I$$
and $[G_{i-1}:G_i]=p_i$.
\rm
\\ \\
See [H, Corollary 10.5.2] or [R, Theorem 5.4.8] for a proof. 
Strangely enough, this theorem is not mentioned in this form 
in the long survey article on supersolvable groups in [Wei]. 
\par
We point out two immediate consequences that will also be used 
frequently in this paper. If $G$ is supersolvable, the Sylow 
$p$-subgroup for the biggest prime $p$ is normal. Also 
$G\cong N\rtimes P$ where $P$ is the Sylow $p$-subgroup for the 
smallest prime that divides $|G|$.
\par
Another fact that we will frequently use is
\\ \\
{\bf Theorem 3.2.} [H, Theorem 10.5.4], [R, Theorem 5.4.10] \it 
If the group $G$ is supersolvable, then its commutator group $G'$ 
is nilpotent. 
\rm
\\ \\
In this paper we not only have to deal with groups which we 
assume to be supersolvable. For some of the results we also 
have to construct groups of a certain form and show that they 
are supersolvable. A very useful tool for this is the following
criterion, which is not completely obvious. 
\\ \\
{\bf Theorem 3.3.} \it
Let $p<q$ be primes.
\begin{itemize}
\item[(a)] If $|G|=p\cdot q^n$ and  $p|(q-1)$, then $G$ is supersolvable.
\item[(b)] If $|G|=p^2 q^n$ and  $p^2 |(q-1)$, then $G$ is supersolvable.
\end{itemize}
\rm

\noindent
{\bf Proof.} \rm 
See [Wei p.6, Corollary 1.10] and the Exercise after it. The key point 
is that these groups have a normal Sylow $q$-subgroup $Q$ and that 
$G/Q$ is abelian of exponent dividing $q-1$. Such groups, which in 
[Wei] are called strictly $q$-closed, are shown to be supersolvable 
in [Wei, p.5 Theorem 1.9].
\hfill$\Box$
\\ \\
A group $G$ is called {\bf metacyclic} if it has a normal cyclic 
subgroup $N$ such that $G/N$ is also cyclic. Note that [H, p.146]
uses a more restrictive definition for metacyclic than we (and most 
sources) do, namely that $G'$ and $G/G'$ are cyclic.
\\ \\
{\bf Lemma 3.4.} \it
Let $G$ be a metacyclic group and $p$ the smallest prime that 
divides $|G|$. Then $G$ has a normal cyclic subgroup $N$ such
that $G/N$ is cyclic and $p$ divides $|G/N|$.
\rm
\\ \\
{\bf Proof.} \rm 
If not, then $G\triangleright N\cong C_{p^e m}$ and $G/N\cong C_n$ 
such that $m$ and $n$ are only divisible by primes that are bigger 
than $p$. In that case $C_m \triangleleft G$ and $G/C_m =H$ has 
a normal subgroup $C_{p^e}$ with quotient $C_n$. By the 
Schur-Zassenhaus Theorem [R, Theorem 9.1.2] we have 
$H\cong C_{p^e} \rtimes C_n$. Moreover, 
$|Aut(C_{p^e})|=p^{e-1}(p-1)$; so $C_n$ can only act trivially 
on $C_{p^e}$. Thus $H\cong C_{p^e n}$.
\hfill$\Box$
\\ \\
A {\bf Z-group} is a finite group whose Sylow subgroups are 
all cyclic. Such groups become important if one wants to bound
the exponent of $G\subseteq Aut(X)$. Compare Theorem 5.8. The 
following relation with other types of groups is not completely 
obvious.
\\ \\
{\bf Theorem 3.5. (Zassenhaus)} [H, Theorem 9.4.3], 
[R, Theorem 10.1.10] \it
A $Z$-group that is not cyclic can be written as a semidirect product
$$C_m\rtimes C_n$$
where $(m,n)=1$ and $m$ is odd. In particular, such a group is 
split metacyclic.
\rm
\\

\subsection*{4. Bounds involving the smallest prime that divides $|G|$}

\noindent
Implicitly the bound from Theorem 2.3 (c) also occurs somewhere
else in the literature for a different type of groups, namely:
\\ \\
{\bf Example 4.1.} Let $p$, $q$ be primes with $p\geq 5$ and 
$p|(q-1)$. The paper [W1] determines every genus on which the
(up to isomorphism unique) non-abelian group $C_q \rtimes C_p$ 
can act. By [W1, Corollary 4.2], the smallest such genus is 
$\mu=1+q\frac{p-3}{2}$. 
\par
So if $C_q \rtimes C_p\subseteq Aut(X)$ we always have 
$|C_q \rtimes C_p|\leq\frac{2p}{p-3}(g-1)$, 
with equality if and only if $g=\mu$.
\par
Since by Dirichlet's Theorem on primes in arithmetic progressions
for fixed $p$ there are infinitely many primes $q$ with $p|(q-1)$,
this gives infinitely many $G$ that are not nilpotent with 
$|G|=\frac{2p}{p-3}(g-1)$.
\\ \\
We recall that $\GGG(p)$ denotes the class of all finite groups 
whose orders are not divisible by any prime smaller than $p$. 
Inspired by the coincidence of the bounds in Theorem 2.3 (c) 
and Example 4.1 we first tried to prove that the same bound 
holds for supersolvable $G$ in $\GGG(p)$ by an inductive process 
using Zappa's Theorem. Then we realized that it actually holds 
for {\it all} $G$ in $\GGG(p)$.
\\ \\
{\bf Theorem 4.2.} \it
Fix a prime $p\geq 5$. Let $X$ be a compact Riemann surface 
of genus $g\geq 2$, and let $G$ be a subgroup of $Aut(X)$ such 
that $|G|$ is not divisible by any primes that are smaller than 
$p$. Then
$$|G|\leq\frac{2p}{p-3}(g-1).$$
Moreover, if $|G|=\frac{2p}{p-3}(g-1)$, then $G$ is a quotient
of $\Gamma(0;p,p,p)$. 
\rm
\\ \\
{\bf Proof.} \rm 
If $G$ is covered by a Fuchsian group 
$\Gamma(h;m_1 ,m_2 ,\ldots, m_r)$, then from Theorem 2.1 
we see that for $|G|\geq\frac{2p}{p-3}(g-1)$ we need 
$2h-2+\sum_{i=1}^r (1-\frac{1}{m_i})\leq 1-\frac{3}{p}$.
Moreover, the periods $m_i$ must divide the group order,
so $m_i \geq p$. This implies $h=0$, and since for $h=0$
we have $r\geq 3$, we see that $\Gamma(0;p,p,p)$ is the 
only possibility. 
\hfill$\Box$
\\ \\
{\bf Remark 4.3.} 
Theorem 2.3 (c) shows that there are infinitely many values of 
$g$ for which the bound in Theorem 4.2 is attained, even if we 
restrict to certain classes of groups, for example $p$-groups, 
nilpotent, or supersolvable.
\par
And Example 4.1 provides infinitely many examples that reach 
the bound if we restrict to groups that are of square-free order 
or $Z$-groups or metacyclic or metabelian. 
\par
On the other hand, $\Gamma(0;p,p,p)$ shows that the only possible 
abelian groups that reach the bound are $C_p$ and $C_p \times C_p$.
\rm
\\ \\
If a group $G$ in Theorem 4.2 reaches the bound 
$|G|=\frac{2p}{p-3}(g-1)$, then $|G|$ obviously must be divisible 
by $p$. Now we formulate the analogues of Theorem 2.2 (b) and (c) 
for such groups.
\\ \\
{\bf Proposition 4.4.} \it
Fix a prime $p\geq 5$. Let $X$ be a compact Riemann surface 
of genus $g\geq 2$ and $G\subseteq Aut(X)$ such that $G$ is 
in $\GGG(p)$ and reaches the bound from Theorem 4.2.
\begin{itemize}
\item[(a)] If $G$ is supersolvable, then $|G|$ is only divisible 
by $p$ and by primes $q$ with $q\equiv 1\ mod\ p$.
\item[(b)] If $G$ is nilpotent, then $G$ must be a $p$-group.
\end{itemize}
\rm

\noindent
{\bf Proof.} \rm 
(a) We know that $H=G/G'$ is isomorphic to $C_p$ or $C_p \times C_p$.
If $G$ is supersolvable, then $G'$ is nilpotent by Theorem 3.2. Fix 
a prime $q\neq p$ that divides $|G'|$ and let $N$ be the product of 
all Sylow subgroups of $G'$ for primes different from $q$. Then $N$ 
is characteristic in $G'$ and hence normal in $G$. So if $Q$ is 
a Sylow $q$-subgroup of $G$, we get a quotient $G/N$ which is 
isomorphic to $Q \rtimes H$ and supersolvable (as a quotient of 
a supersolvable group). 
By Zappa's Theorem, this quotient has a normal subgroup which has 
index $q$ in $Q$. So we see that $G$, and hence $\Gamma(0;p,p,p)$ 
has a quotient of the form $C_q \rtimes H$. If $q$ is not congruent 
to $1$ modulo $p$, this group is abelian, contradicting $H=G/G'$.
\par
(b) If $G$ is nilpotent, then any Sylow $q$-subgroup of $G$ is
a quotient of $G$ and hence of $\Gamma(0;p,p,p)$, forcing $q=p$.
\hfill$\Box$
\\ \\
Proposition 4.4 (b) is an anlogue of Theorem 2.2 (c). 
For the corresponding statement for nilpotent groups of odd 
order see Theorem 5.2. Example 4.1 shows that the statement
'equality only possible for $p$-groups' does not generalize 
to supersolvable $G$.
\par
We construct more examples of $G$ that are metabelian {\it and} 
supersolvable and attain the bound from Theorem 4.2.
\\ \\
{\bf Example 4.5.} 
Let $p\geq 5$ be a prime. Then there exists a Riemann surface of
genus $h=\frac{p-1}{2}$ with $C_p \subseteq Aut(X)$. Fix any prime
$q>p$. 
By Macbeath's construction, for every $e>0$ there is a smooth Galois 
cover $X$ of $Y$ whose Galois group is a direct product of $2h$ copies 
of $C_{q^e}$. So the genus of $X$ is $g=\frac{p-3}{2}q^{2he}+1$. The 
automorphism group of $Y$ lifts to $X$. So $Aut(X)$ has a subgroup 
$G\cong (C_{q^e})^{2h} \rtimes C_p$. Thus $G$ is metabelian 
by construction and attains the bound from Theorem 4.2.
\par
By Dirichlet's Theorem on primes in arithmetic progressions there 
are infinitely many primes $q$ with $q\equiv 1\ mod\ p$. If we take 
one of those, $G$ moreover is supersolvable by Theorem 3.3 (a).
\rm
\\ \\
{\bf Remark 4.6.} 
From Theorems 2.2 (b) and Examples 4.1 and 4.5 one might get 
the impression that if $G\subseteq Aut(X) $ is a supersolvable, 
non-nilpotent group of maximal possible order, then the smallest 
prime can divide $|G|$ only once. Although this is true if the 
smallest prime is $2$ or $3$ (see Theorem 5.1 below), there are 
infinitely many counterexamples for every $p\geq 5$.
\par
We start with the group of order $p^2$ acting on a surface $Y$ of 
genus $h=\frac{p(p-3)}{2}+1$ (see Theorem 2.3 (c)). Being a quotient 
of $\Gamma(0;p,p,p)$, it can only be $C_p \times C_p$. Then, by the
same construction as in Example 4.5 for every prime $q>p$ and every
$e>0$ we obtain a metabelian group 
$G\cong (C_{q^e})^{2h} \rtimes (C_p \times C_p)$
acting on a surface of genus $g=\frac{p(p-3)}{2}q^{2he}+1$. 
\par
If $q$ is one of the infinitely many primes with 
$q\equiv 1\ mod\ p^2$, then moreover $G$ is supersolvable 
by Theorem 3.3 (b).
\rm
\\ \\
{\bf Theorem 4.7.} \it
Let $p$ be an odd prime. If $G\subseteq Aut(X)$ is cyclic 
and in $\GGG(p)$, then 
$$|G|\leq \frac{2p}{p-1}g+p.$$
This bound is attained for every $g=\frac{p-1}{2}(m-1)$
with a group $C_{pm}$, provided all prime divisors of
$m>1$ are bigger than $p$.
\rm
\\ \\
{\bf Proof.} \rm 
Let $G$ be a cyclic group of order $N$ and $p_1$ the smallest 
prime divisor of $N$. Then by [Ha, Theorem 6] the minimum genus
on which $G$ can act is $\frac{p_1 -1}{2}\cdot\frac{N}{p_1}$
if $N$ is prime or $p_1^2$ divides $N$, and 
$\frac{p_1 -1}{2}\cdot(\frac{N}{p_1}-1)$ otherwise. 
This proves our theorem provided $p$ divides $|G|$.
\par
Now assume that the smallest prime divisor of $|G|$ is $q>p$.
If $|G|$ is bigger than the bound in our theorem, then 
necessarily $\frac{2p}{p-1}g+p<\frac{2q}{q-1}g+q$. An easy
algebraic manipulation shows that this inequality is equivalent
to $g<\frac{(p-1)(q-1)}{2}$. Plugging this into the bigger bound
yields $|G|<pq$, leaving only the possibility $|G|=q$, for which 
the stronger $|G|\leq 2g+1$ holds.
\hfill$\Box$
\\ \\
{\bf Theorem 4.8.} \it
Let $p$ be an odd prime. If $G\subseteq Aut(X)$ is abelian 
and in $\GGG(p)$, then 
$$|G|\leq \frac{2p}{p-1}g+2p,$$
except for the groups $C_q \times C_q$ where $q$ is a prime
with $p<q<2p$ acting on genus $g=\frac{(q-1)(q-2)}{2}$.
The bound can be reached for every $g=\frac{p-1}{2}(m-2)$
for which $p$ is the smallest prime divisor of $m$ with 
a group $C_p \times C_m$.
\rm
\\ \\
{\bf Proof.} \rm 
Let $A\cong C_{m_1}\times C_{m_2}$ with $m_1 |m_2$. Then by 
[Ml, Theorem 4] the minimum genus $g^*$ on which $A$ can 
act is given by 
$\frac{2(g^*-1)}{|A|}=1-\frac{1}{m_1}-\frac{2}{m_2}$.
This transforms to $|A|\leq \frac{2m_1}{m_1 -1}g+2m_1$,
and it also shows that for $m_1 =p$ the bound in 
Theorem 4.8 can be reached. 
\par
That bound obviously surpasses the bound from 
Theorem 4.7 for cyclic groups in $\GGG(p)$. On the other 
hand, for abelian groups of odd order and rank $r>2$ by
[BMT, Theorem 2] we have $|A|\leq 2g+6$, except for 
$C_3 \times C_3 \times C_3$, which can act on genus $10$. 
\par
So if an abelian group $A$ in $\GGG(p)$ tops the bound
in the Theorem 4.8, it must be of the form 
$A\cong C_{m_1}\times C_{m_2}$ with $m_1 |m_2$ and $m_1 >p$.
Moreover, we then must have
$\frac{2p}{p-1}g+2p < \frac{2m_1}{m_1 -1}g+2m_1$, which
is equivalent to $g<(p-1)(m_1 -1)$. Then the bigger bound
gives $|A|<2pm_1$. Consequently, $A$ can only be 
$C_q \times C_q$ with a prime $q$ between $p$ and $2p$.
By [Ml, Theorem 4] these groups can indeed act on genus
$g^*=\frac{(q-1)(q-2)}{2}$ and hence surpass the bound in
the theorem. The next bigger genus on which they can act
is $g^* +q$, which is too big to compete with the bound 
from the theorem. Here we are using that 
$\frac{|A|}{exp(A)}=q$ must divide $2(g-1)$, 
a well-known fact that was reproved and elaborated 
on in [Sch1, Section 6].
\hfill$\Box$
\\

\subsection*{5. Groups of odd order}

\noindent 
In this section we refine Theorem 2.4 by working out the
analogue of Table 1 for groups of odd order.
\\ \\
{\bf Theorem 5.1.} \it 
Let $X$ be a compact Riemann surface of genus $g\geq 2$, and let $G$
be a subgroup of $Aut(X)$. If $G$ is supersolvable and of odd order, 
then $|G|\leq\frac{21}{2}(g-1)$.
\par 
Moreover, a supersolvable group $G$ of odd order that reaches 
this bound necessarily has $|G|=3\cdot 7^n$. 
\par
Conversely, for every $n\geq 1$ there exists a Riemann surface
$X$ of genus $g=2\cdot 7^{n-1}+1$ such that $Aut(X)$ contains 
a supersolvable subgroup $G$ of order $|G|=3\cdot 7^n$.
\rm
\\ \\
{\bf Proof.} \rm 
To show $|G|\leq\frac{21}{2}(g-1)$ by Table 3 we only have to exclude 
the possibility $|G|=15(g-1)$. So suppose $|G|=15(g-1)$, i.e., that
$G$ is a quotient of $\Gamma(0;3,3,5)$. Then by Zappa's Theorem
$G$ has a quotient of order $9$ (if $9$ divides $|G|$) or a quotient
of order $15$ (if $9$ does not divide $|G|$). Either one contradicts
$G/G'\cong C_3$. 
\par
Now let $G\subseteq Aut(X)$ be a supersolvable group of odd order
$\frac{21}{2}(g-1)$. Then $G$ is a quotient of $\Gamma(0;3,3,7)$.
By Table 3 we have $G/G'\cong C_3$ and $G'$ is a quotient of 
$\Gamma(0;7,7,7)$. Moreover, $G'$, as the commutator group of 
a supersolvable group, must be nilpotent by Theorem 3.2. Thus
every Sylow $p$-subgroup of $G'$ is a quotient of $G'$ and
hence of $\Gamma(0;7,7,7)$. But obviously this group has no 
quotients of order prime to $7$. So $|G'|=7^n$, and hence 
$|G|=3\cdot 7^n$.
\par
There is a non-abelian group of order $21$ acting on a Riemann
surface $X$ of genus $3$. More precisely, $X$ is the Klein quartic
and $G$ is the normalizer of a Sylow $7$-subgroup in 
$Aut(X)\cong PSL_2(\FF_7)$.
\par 
By [W2, Theorem 4.2] for every $\mu\geq 0$ there is a group $G$ 
of order $3\cdot 7^{30\mu +2}$ acting on a Riemann surface $X$ of 
genus $g=2\cdot 7^{30\mu +1}+1$. These $G$ are supersolvable by 
Theorem 3.3. To fill in the groups of order $3\cdot 7^n$ with 
the missing $n$, we start with such a group $G$ for a big enough 
$n=30\mu +2$ and inductively divide the group order by $7$ as 
follows. 
By Zappa's Theorem, $G$ has a normal subgroup $P\cong C_7$. Then
$G/P$ acts on the Riemann surface $X/P$. If $h$ denotes the genus 
of $X/P$, then $h-1\leq \frac{1}{7}(g-1)$ by the Hurwitz formula.
Because of the bound $|G/P|\leq\frac{21}{2}(h-1)$, which we have 
already established in general, we actually must have 
$h-1=\frac{1}{7}(g-1)$, except if it unfortunately should happen 
that $h\leq 1$. But in that case $G'/P$, a group of order $7^{n-1}$ 
acting on a surface of genus $0$ or $1$, would have to be abelian,
which for $n\geq 4$ contradicts the fact that $C_7 \times C_7$ is
the biggest abelian quotient of $\Gamma(0;7,7,7)$.
\hfill$\Box$
\\ \\
See also Proposition 7.3 and Theorem 8.2 for further 
characterizations of the groups in Theorem 5.1.
\par
The following analogue of Theorem 2.2 (c) is at least implicitly
already in [Z1] or [Z2].
\\ \\
{\bf Theorem 5.2.} \it
Let $X$ be a compact Riemann surface of genus $g\geq 2$, and 
let $G$ be a subgroup of $Aut(X)$. If $G$ is nilpotent and of 
odd order, then $|G|\leq 9(g-1)$. Moreover, any nilpotent group 
of odd order that reaches this bound must be a $3$-group.
\rm
\\ \\
{\bf Proof.} \rm 
If $G$ is nilpotent, then every Sylow $p$-subgroup of $G$ is 
a quotient of $G$. But the triangle groups $\Gamma(0;3,3,n)$ 
have no nontrivial quotients of order prime to $3$. This rules 
out the possibilities $|G|=15(g-1)$ and $|G|=\frac{21}{2}(g-1)$. 
The next one is $|G|=9(g-1)$ with $G$ being a quotient of 
$\Gamma(0;3,3,9)$. By the same argument then $|G|$ is not 
divisible by any primes $p\neq 3$.
\hfill$\Box$
\\ \\
{\bf Theorem 5.3.} \it
If $G\subseteq Aut(X)$ is a metabelian group of odd order, 
then, apart from $3$ exceptions we have $|G|\leq 9(g-1)$.
\par
The exceptions are the non-abelian group of order $21$ acting
on a surface of genus $3$, a group $(C_7 \times C_7)\rtimes C_3$ 
for $g=15$, and a group $(C_5 \times C_5)\rtimes C_3$ for $g=6$.
\par
Moreover, the bound $9(g-1)$ is attained infinitely often.
More precisely, if $m\in \NN$ is only divisible by primes that
are congruent to $1$ modulo $3$, then there exists a compact 
Riemann surface $X$ of genus $g=9m+1$ such that $Aut(X)$ 
contains a subgroup $G$ of order $9(g-1)$ that is supersolvable
and metabelian.
\rm
\\ \\
{\bf Proof.} \rm 
If $|G|$ is bigger than $9(g-1)$, then $G$ is covered by 
$\Gamma(0;3,3,p)$ with $p=7$ or $p=5$. So $G/G'\cong C_3$
and $G'$ is a quotient of $\Gamma(0;p,p,p)$. Since $G'$ 
is abelian, this leaves only the possibilities $G'\cong C_p$
or $G'\cong C_p \times C_p$. Moreover $G'\cong C_5$ is not 
possible, as then $G$ would be cyclic of order $15$ on 
a surface of genus $2$.
\par
To prove the last statement, we note that by 
Theorem 2.2 (b) for the specified genus $g=9m+1$ 
there exists a surface $X$ with supersolvable 
$H\subseteq Aut(X)$ of order $18(g-1)$. 
We take $G$ to be the index $2$ subgroup in $H$. 
Then $G$ has the right order and is supersolvable. So there 
only remains to show that $G$ is metabelian. 
\par
Note that $G$ is a quotient of $\Gamma=\Gamma(0;3,3,9)$.
We have $\Gamma/\Gamma'\cong C_3 \times C_3$ and 
$\Gamma'=\Gamma(1;3,3,3)$. By Zappa's Theorem $G$ must have 
a quotient of order $9$. So $G/G'\cong C_3 \times C_3$ and
$G'$ is a quotient of $\Gamma(1;3,3,3)$. Moreover, as $G$ is 
supersolvable, $G'$ must be nilpotent by Theorem 3.2. 
So $G'$ is a direct product of its Sylow subgroups. The Sylow
$3$-subgroup $N$ is of order $9$ and in particular abelian.
Actually, $G'/N$, being a quotient of $\Gamma(1;3,3,3)$ of 
order prime to $3$, must also be abelian. Thus $G'$ is abelian,
which finishes the proof.
\hfill$\Box$
\\ \\
For more information on the three exceptions see [MZ2].
\\ \\
{\bf Theorem 5.4.} \it 
\begin{itemize}
\item[(a)] If $G\subseteq Aut(X)$ is a cyclic group of odd 
order, then 
$$|G|\leq 3g+3.$$ 
This bound is attained for every $g\equiv\ 0$ or $4\ mod\ 6$.
\item[(b)] If $G\subseteq Aut(X)$ is an abelian group of odd 
order, then with the exception of $C_5 \times C_5$ acting on
a Riemann surface of genus $6$ we have  
$$|G|\leq 3g+6.$$ 
This bound can be attained for every $g=6k+1$ where 
$k=1,2,\ldots$ with a group $G\cong C_3 \times C_{6k+3}$.
\end{itemize}
\rm

\noindent
{\bf Proof.} \rm 
(a) is is a special case of Theorem 4.7. Alternatively, 
it follows directly from [N].
\par
(b) is is a special case of Theorem 4.8. 
\hfill$\Box$
\\ \\
Next we prove the odd order analogues of the very recent 
results in [Sch2] concerning several types of metacyclic 
subgroups of $Aut(X)$.
\\ \\
{\bf Theorem 5.5.} \it
If $G\subseteq Aut(X)$ is a metacyclic group of odd order, 
then
$$|G|\leq 6g+3.$$
This bound can be attained if and only if $2g+1$ is neither
divisible by $9$ nor by any prime $p$ with $p\equiv 2\ mod\ 3$.
The corresponding $G$ then is of the form $C_{2g+1}\rtimes C_3$.
\rm
\\ \\
{\bf Proof.} \rm 
If $G$ is in $\GGG(5)$, then $|G|\leq 5(g-1)$ by Theorem 4.2.
So we can assume that $3$ divides $|G|$. By Lemma 3.4 then $G$
has a normal cyclic subgroup $N\triangleleft G$ such that 
$G/N\cong C_n$ with $3|n$. Actually, we must have $n=3$. Bigger 
cyclic quotients are not possible, because then the periods 
would make $|G|$ too small. For example, $\Gamma(0;3,5,15)$ 
gives $|G|=5(g-1)$. Alternatively, for $n>9$ one could invoke 
[MZ1, Theorem 1].
\par
This also immediately shows that $|N|$ cannot be divisible
by any prime $p$ with $p\equiv 2\ mod\ 3$, for then we would 
have $G/H\cong C_{3p}$ where $H$ is the unique subgroup of
index $p$ in $N$. 
\par
So if $|G|\geq 6g+3$, then $G$ has a normal cyclic subgroup 
$N\cong C_m$ with $m\geq 2g+1$. Thus $N$ is what is called 
a quasilarge abelian group (of automorphisms) in [PR]. 
By [PR, Proposition 4.1] the quotient surface $X/N$ must 
have genus $0$. More precisely, $N$ can only be an instance
of the last case (3 critical values) in [PR, Table 2]. 
So $N$ is covered by some $\Gamma(0;m_1 ,m_2 , m_3 )$.
\par
Now by [BC, Theorem 4.2] the fact that $N$ extends normally 
to $G\subseteq Aut(X)$ with $[G:N]=3$ is only possible if 
$m_1 =m_2 =m_3 =t\geq 4$ (see Case N6 in [BC, p.576]). In 
terms of [PR, Table 2] this means $\alpha=1$, 
$\delta_i =1$ and $\beta=t=|N|$. Hence by 
[PR, Theorem 4.2] we have $g=1+\frac{1}{2}(t-3)$, 
i.e. $t=2g+1$. This establishes the bound $|G|\leq 6g+3$.
\par
Now let $z_1$, $z_2$, $z_3\in N$ be the images of the generators
$x_1$, $x_2$, $x_2$ of $\Gamma(0;t,t,t)$ where $t=2g+1$, and let
$u$ be a generator of $G/N\cong C_3$. Then $uz_3u^{-1}=z_3^b$ for 
some integer $b$ that is relatively prime to $t$. By [BC, Case N6]
the necessary and sufficient condition for the existence of $G$
is that conjugation with $u$ permutes the $z_i$ cyclically, i.e.
that $z_1 =z_3^b$ and $z_2 =z_3^{b^2}$. Because of $z_1 z_2 z_3=1$ 
this forces $1+b+b^2$ to be divisible by $t$. From this we see 
that $t$ cannot be divisible by any prime $p\equiv 2\ mod\ 3$
(as we have already seen earlier) and that $9$ cannot divide $t$.
\par
Conversely, if $t=2g+1$ is a product of primes that are congruent 
to $1$ modulo $3$, or $3$ times such a number, then by the Chinese
Remainder Theorem there exists an integer $b$, relatively prime 
to $t$, such that $t$ divides $1+b+b^2$. Then by [BC, Case N6] the 
desired $G$ exists.
Note that in the latter case $b\equiv 1\ mod\ 3$ and the Sylow 
$3$-subgroup of $G$ will then be $C_3 \times C_3$.
\par
Finally, the fact that $G$ cannot have a quotient $C_9$ guarantees
that $N$ has a complement, and hence $G$ is a semidirect product.
\hfill$\Box$
\\ \\
{\bf Remark 5.6.} 
If $p$ is a prime that is congruent to $1$ modulo $3$, then 
by [W1, Corollary 4.2] the smallest genus on which the (up 
to isomorphism unique) non-abelian group $C_p \rtimes C_3$ 
can act is $g=\frac{p-1}{2}$. So the group order is $6g+3$.
This was the initial inspiration for Theorem 5.5.
\par
However, the more general construction with groups 
$C_{2g+1}\rtimes C_3$ in [W2, Section 6] is not completely correct.
For example, $l=4$, $d=1$ gives $g=17$. So the group, which can 
only be $C_5 \times (C_7 \rtimes C_3)$, also has elements of order 
$15$, contrary to what is claimed in [W2, p.219]. Indeed, by 
Theorem 5.5 this group cannot act on a surface of genus $17$. 
\par
Coincidentally, [KS, Table 1], which allows to write down
algebraic equations for examples as in Theorem 5.5, contains
a similar inaccuracy. The conditions in Case C.1, namely that 
$n$ in $C_n\rtimes C_3$ must be bigger than $7$, odd, and 
divisible by a prime that is congruent to $1$ modulo $3$, 
are at least misleading. If $1<b<n$ and $1+b+b^2$ is divisible 
by $n$ (which is the set-up in [KS]), this implies exactly the 
conditions from Theorem 5.5, namely $9\tn n$ and $p\tn n$ for 
every prime $p$ with $p\equiv 2\ mod\ 3$.
\rm
\\ \\
We also point out that the approach in the proof of 
Theorem 5.5 can also be used to prove $|G|\leq 12(g-1)$
for general metacyclic groups independent of the line of
argument in [Sch2]. 
\par
First one uses [MZ1, Theorem 1] to show $|G|\leq 12(g-1)$
if $G$ has a cyclic quotient $C_q$ with $q\geq 7$. (Here 
$q$ is an integer, not necessarily a prime.) So if $G$ is 
metacyclic and $|G|>12(g-1)$, then $G$ has a normal cyclic 
subgroup $N$ of order bigger than $2(g-1)$. Thus $N$ is 
quasilarge and must therefore show up in [PR, Table 2]. 
Then one compares with [BC]. Admittedly, some care with 
the details is required.
\par 
Theorem 5.5, when combined with Theorem 3.5, immediately 
yields more results.
\\ \\
{\bf Corollary 5.7.} \it 
Let $g\geq 2$ and $G\subseteq Aut(X)$.
\begin{itemize}
\item[(a)] If $G$ is a $Z$-group of odd order, then 
$|G|\leq 6g+3$. This bound can be attained if and only 
if all prime divisors of $2g+1$ are congruent to $1$ 
modulo $3$.
\item[(b)] If  $|G|$ is square-free and odd, then 
$|G|\leq 6g+3$. This bound can be attained if and 
only if $2g+1$ is a product of distinct primes that 
all are congruent to $1$ modulo $3$.
\end{itemize}
In both cases, if the bound is attained, then
$G\cong C_{2g+1} \rtimes C_3$.
\rm
\\ \\
Using the easy fact that the exponent $exp(G)$ of a finite 
group $G$ equals $|G|$ if and only if $G$ is a $Z$-group, 
we also obtain the following result. 
\\ \\
{\bf Theorem 5.8.} \it
If $G\subseteq Aut(X)$ has odd order, then 
$$exp(G)\leq 6g+3.$$
This bound can be attained if and only if all 
prime divisors of $2g+1$ are congruent to $1$ 
modulo $3$ (with $G\cong C_{2g+1}\rtimes C_3$).
\rm
\\ \\
{\bf Proof.} \rm 
If $G$ is a $Z$-group, this follows from 
Corollary 5.7 (a). If $G$ is not a $Z$-group, we have 
$exp(G)\leq \frac{1}{3}|G|\leq \frac{1}{3}15(g-1)$
by Theorem 2.4.
\hfill$\Box$
\\ \\
The analogue of Theorem 5.8 for groups of not necessarily odd 
order is [Sch1, Theorem 4.4]. The corresponding problem for 
solvable groups of not necessarily odd order is not completely 
solved. See [Sch2, Proposition 6.2 and Remark 6.3].
\\

\subsection*{6. Groups of order $p^m q^n$}

\noindent
Several times we have seen in the preceding sections that one can 
often find $(p,q)$-groups that reach the sharp bound for a certain 
type of groups $G\subseteq Aut(X)$. Motivation enough for us to 
investigate this type of groups in their own right. Note that
$(p,q)$-groups are always solvable, a famous theorem by Burnside 
[H, Theorem 9.3.2] or [R, Theorem 8.5.3].
\par
In $|G|=p^m q^n$ we always choose the letter $p$ for the smaller 
one of the two primes $p$ and $q$. Moreover, to avoid certain
trivialities, we implicitly assume that $m$ and $n$ are both 
positive. In accordance with Sections 4 and 5 we again refine 
the statements for odd group order and in terms of the smallest 
prime $p$.
\\ \\
{\bf Theorem 6.1.} \it 
Let $X$ be a compact Riemann surface of genus $g\geq 2$, and
let $G\subseteq Aut(X)$ be a subgroup. If $G$ is a $(p,q)$-group,
then $|G|\leq 48(g-1)$. Any $(p,q)$-group that reaches this bound 
necessarily is a $(2,3)$-group. 
\par
Conversely, for every $g=2^{6\mu+1}3^{6\nu}+1$ and every 
$g=2^{4\mu}3^{4\nu}+1$ where $\mu,\nu\geq 0$ there exists 
a Riemann surface $X$ of genus $g$ such that $Aut(X)$ 
is a $(2,3)$-group of order $48(g-1)$.
\rm
\\ \\
{\bf Proof.} \rm 
Since $(p,q)$-groups are solvable, the first claim follows from 
part (a) of Theorem 2.2. The second claim is clear, because then 
$48$ divides $|G|$. For the third claim we use the parametrizations
from the proof of Theorem 2.2 (a) with $n=2^{\mu}3^{\nu}$. 
\hfill$\Box$
\\ \\
{\bf Theorem 6.2.} \it 
Let $X$ be a compact Riemann surface of genus $g\geq 2$, and
let $G\subseteq Aut(X)$ be a subgroup. If $G$ is a $(p,q)$-group
of odd order, then $|G|\leq 15(g-1)$. Any $(p,q)$-group that reaches 
this bound necessarily is a $(3,5)$-group. 
\par
Conversely, for every $g=3^{12\mu}5^{12\nu +1}+1$ with $\mu,\nu\geq 0$ 
there exists a Riemann surface $X$ of genus $g$ such that $Aut(X)$ 
contains a $(3,5)$-group $G$ of order $15(g-1)$.
\rm
\\ \\
{\bf Proof.} \rm 
The bound and the existence follow from 
Theorem 2.4 (with $m=3^{\mu}5^{\nu}$).
\hfill$\Box$
\\ \\
In Example 4.5 and Remark 4.6 we already showed the following result. 
\\ \\
{\bf Theorem 6.3.} \it
Fix primes $p$, $q$ with $3<p<q$. Then there are infinitely many 
metabelian $(p,q)$-groups $G\subseteq Aut(X)$ that reach the bound 
$|G|=\frac{2p}{p-3}(g-1)$ from Theorem 4.2. 
\rm
\\ \\
Applying Macbeath's construction to the $p$-groups in Theorem 2.3 (c), 
one can for every $n>0$ (and every $q$) obtain $(p,q)$-groups $G$ with 
$ord_p |G|=n$ that attain the bound. But then in general $G$ will be 
neither metabelian nor supersolvable.
\par
Finally, we address the problem what happens if we fix $p\leq 3$ and 
$q$. Actually, $(2,3)$-groups were already discussed in Theorem 6.1.
\\ \\
{\bf Theorem 6.4.} \it 
Fix an odd prime $q$. There are infinitely many metabelian 
$(2,q)$-groups $G\subseteq Aut(X)$ with $|G|=16(g-1)$.
For $q\geq 11$ this is also the upper bound that can be 
achieved by a $(2,q)$-group.
\par
For $(2,5)$-groups (resp. $(2,7)$-groups) $G\subseteq Aut(X)$ 
we have $|G|\leq 40(g-1)$ (resp. $|G|\leq \frac{56}{3}(g-1)$),
and this bound is reached infinitely often.
\rm
\\ \\
{\bf Proof.} \rm 
Setting $\beta=2$ and $k=2^{\mu}q^{\nu}$ in [G3, Theorem 1.2] 
we get a metabelian group of order $2^{\mu +5}q^{\nu}$ acting 
on a surface of genus $g=2^{\mu +1}q^{\nu}+1$. 
\par
All orders in Table 2 are divisible by $3$, $5$ or $7$. 
There are a handful more possible orders that are bigger than 
$16(g-1)$, to wit, for $\Gamma(0;2,3,m)$ with $19\leq m\leq 23$, 
but they are all divisible by the wrong primes. 
\par
There are Riemann surfaces $X$ of genus $5$ with $|Aut(X)|=160$;
they are called Humbert curves. By Macbeath's construction, from 
them we can obtain for any $\mu, \nu \geq 0$ a group of order 
$2^{10\mu +5}5^{10\nu +1}$ acting on a Riemann surface of genus 
$g=4\cdot(2^{\mu}5^{\nu})^{10}+1$. By [ChP] these groups cannot 
be metabelian.
\par
Finally to $|G|=\frac{56}{3}(g-1)$. Let $\Gamma=\Gamma(0;2,4,7)$.
Then $\Gamma/\Gamma'\cong C_2$ and $\Gamma'=\Gamma(0;7,7,2)$.
So $\Gamma'/\Gamma''\cong C_7$ and 
$\Gamma''=\Gamma(0;2,2,2,2,2,2,2)$. Finally, 
$\Gamma''/\Gamma'''\cong (C_2)^6$ and $\Gamma'''=\Gamma(49;-)$.
So there exists a Riemann surface $X$ of genus $g=49$ with 
a $G\subseteq Aut(X)$ with $|G|=7\cdot 2^7=\frac{56}{3}(g-1)$.
To that one we can apply Macbeath's construction and get infinitely
many more. However, the genus grows rapidly. The next one constructed
by that method has genus $48\cdot 2^{98}+1$.
\hfill$\Box$
\\ \\
{\bf Theorem 6.5.} \it 
The bound on the order of a $(3,q)$-group $G$ in $Aut(X)$ is 
$15(g-1)$ if $q=5$, $\frac{21}{2}(g-1)$ if $q=7$, and $9(g-1)$
if $q\geq 11$. For any fixed $q$ this bound is attained infinitely 
often.
\rm
\\ \\
{\bf Proof.} \rm 
This follows from the fact that $15(g-1)$, $\frac{21}{2}(g-1)$ 
and $9(g-1)$ are the three smallest possible odd orders together 
with Theorems 2.4 and 5.1 respectively applying Macbeath's 
construction to the $3$-groups from Theorem 2.3 (b).
\hfill$\Box$
\\ \\
If $q$ is congruent to $1$ modulo $3$, then Theorem 5.3 gives 
the much stronger result that for every $n\geq 1$ there exists
a metabelian, supersolvable group $G\subseteq Aut(X)$ of order
$3^4 q^n =9(g-1)$.
\\ \\
{\bf Remark 6.6.} 
If we only fix the bigger prime $q$, then we obtain from the 
previous results that the biggest order of a $(p,q)$-group 
$G\subseteq Aut(X)$ is always obtained with $(2,q)$-groups
as given earlier in this section.
\rm
\\ \\

\subsection*{7. Between supersolvable and solvable}

\noindent
The converse of Theorem 3.2 is not true. Groups $G$ with 
nilpotent commutator group $G'$ are of course solvable, 
but not necessarily supersolvable, as can be seen from 
the standard counterexample $A_4$. On the other hand, 
$S_4$ shows that not every solvable $G$ has a nilpotent 
$G'$. So this type of groups lies properly between the 
supersolvable and the solvable ones. More or less by 
definition it also lies strictly between the metabelian 
and the solvable groups.
It turns out that the optimal bound for $|G|$ in terms of 
$g$ for this type of groups also lies strictly between the 
bounds for supersolvable or metabelian and solvable $G$.  
\\ \\
{\bf Theorem 7.1.} \it 
Let $X$ be a compact Riemann surface of genus $g\geq 2$, and
let $G\subseteq Aut(X)$ be a subgroup. If the commutator group
$G'$ is nilpotent, then $|G|\leq 24(g-1)$. 
\par
Conversely, for every $g=2^{6n+1}+1$ and for every 
$g=2^{4n}+1$ there exists an $X$ and $G\subseteq Aut(X)$ with 
$|G|=24(g-1)$ such that $G/G'\cong C_3$ and $G'$ is a $2$-group.
\rm
\\ \\
{\bf Proof.} \rm 
By Table 2 there are only $4$ possible orders bigger than 
$24(g-1)$. Of these, $48(g-1)$, coming from $\Gamma(0;2,3,8)$, 
can only have $G/G'\cong C_2$; so $G'$, of order $24(g-1)$, 
cannot be nilpotent by Theorem 2.2 (c).
The case $|G|=36(g-1)$ corresponds to $\Gamma(0;2,3,9)$. So 
if $G'$ is a proper subgroup, we have $G/G'\cong C_3$ and 
by Table 2 $G'$ is a quotient of $\Gamma(0;2,2,2,3)$.
In particular, $G'$ cannot have a quotient that is a $3$-group.
But since $|G'|=12(g-1)$ is divisible by $3$, this means that 
$G'$ cannot be nilpotent. The same proofs work for $40(g-1)$ 
and $30(g-1)$, respectively.
\par
For the existence proof we start with a solvable group 
$H$ of order $48(g-1)$ where $g=2^{6n+1}+1$ or $g=2^{4n}+1$. 
By the proof of Theorem 2.2 (a) such groups exist for every 
$n\geq 0$. Then $G=H'$ is the desired group. It has order 
$24(g-1)$ and is a quotient of $\Gamma(0;3,3,4)$ 
(see Table 2). So $G/G'\cong C_3$, which for our choices 
ensures that $G'$ is a $2$-group.
\hfill$\Box$
\\ \\
In contrast, groups of order $3\cdot 5^n$ obviously have 
a nilpotent commutator group, and we know from Theorem 2.4
that there are infinitely many among them that reach the bound
$15(g-1)$ for general (that is, solvable) groups of odd order. 
We now show that among the groups of odd order that reach this 
bound they are the only ones with nilpotent commutator group.
\\ \\
{\bf Proposition 7.2.} \it 
Let $G$ be a group of odd order that reaches the bound
$|G|=15(g-1)$ in Theorem 2.4. Then the commutator group
$G'$ is nilpotent if and only if $|G|=3\cdot 5^n$.
\rm
\\ \\
{\bf Proof.} \rm 
If $G$ is of odd order with $|G|=15(g-1)$, then $G$ is a quotient 
of $\Gamma=\Gamma(0;3,3,5)$. So $G/G' \cong C_3$ and $G'$ is a 
quotient of $\Gamma'=\Gamma(0;5,5,5)$. This shows that $G'$ cannot 
have any quotients of order prime to $5$. So $G'$, if nilpotent, 
must be a $5$-group.
\hfill$\Box$
\\ \\
In the same vein we can show that the second statement of 
Theorem 5.1 holds under slightly weaker conditions.
\\ \\
{\bf Proposition 7.3.} \it 
Let $X$ be a compact Riemann surface of genus $g\geq 2$. 
Let $G\subseteq Aut(X)$ be a subgroup of odd order with 
$|G|=\frac{21}{2}(g-1)$. Then the following are equivalent:
\begin{itemize}
\item[(a)] $|G|=3\cdot 7^n$.
\item[(b)] $G$ is supersolvable.
\item[(c)] The commutator group $G'$ is nilpotent.
\end{itemize}
\rm

\noindent
{\bf Proof.} \rm 
By Theorem 3.3, groups of order $3\cdot 7^n$ are supersolvable. 
And supersolvable groups have a nilpotent commutator group by
Theorem 3.2. So we only have to show that (c) implies (a), which 
is exactly the same argument as in the proof of Proposition 7.2.
\hfill$\Box$
\\ \\
See also Theorem 8.2 for yet another characterization of the 
groups in Theorem 5.1 and Proposition 7.3.
\\ \\
Another property of supersolvable groups is that the elements
of odd order form a (normal) subgroup. This is an immediate 
consquence of Zappa's Theorem. The converse is not true. 
For example it is well known that in a group $G$ with 
$|G|\equiv 2\ mod\ 4$ the elements of odd order form a normal 
subgroup of index $2$. But such a group need of course not be 
supersolvable. 
\par
Now let $G$ be a finite group in which the elements of odd 
order form a subgroup $H$. Then, obviously, $H$ must be normal.
Moreover, $H$ is solvable by the Theorem of Feit and Thompson 
[FT]. Since $G/H$, being a $2$-group, is also solvable, we see
that $G$ must be solvable. But of course not every solvable 
group has the property that the elements of odd order form 
a subgroup.
\\ \\
{\bf Theorem 7.4.} \it
Let $X$ be a compact Riemann surface of genus $g\geq 2$.
Let $G\subseteq Aut(X)$ such that the elements of odd order 
in $G$ form a subgroup $N$. Then $|G|\leq 30(g-1)$ and any 
$G$ that attains this bound must necessarily have 
$|G|\equiv 2\ mod\ 4$.
Conversely, there are infinitely many $g$, necessarily
with $g\equiv 6\ mod\ 10$, for which this bound is sharp.
\rm
\\ \\
{\bf Proof.} \rm 
We have $G/N\cong P$ where $P$ is the Sylow $2$-subgroup.
If $4$ divides $|G|$, then $P$ and hence $G$ has an abelian 
quotient of order $4$. This excludes the three biggest 
possibilities from Table 2. 
\par
If $|G|=30(g-1)$, for the same reason we must have
$|G|\equiv 2\ mod\ 4$. In other words, $g-1$ must be odd.
The corresponding triangle group is 
$\Gamma=\Gamma(0;2,3,10)$. Then $\Gamma'=\Gamma(0;3,3,5)$ 
and $\Gamma''=\Gamma(0;5,5,5)$. So $G'/G''\cong C_3$. Since 
$G''/G'''$ cannot also be cyclic by [H, Theorem 9.4.2],
we must have $G''/G'''\cong C_5 \times C_5$, and hence
$5|(g-1)$. So together $g\equiv 6\ mod\ 10$.
\par
By the construction in [Ch, Theorem 3.2] for $m=5$ and 
any odd $n$ we get $|G|=2\cdot 3\cdot 5^2\cdot n^{12}$ 
on a Riemann surface of genus $5\cdot n^{12}+1$.
\hfill$\Box$
\\ \\
{\bf Corollary 7.5.} \it
Let $X$ be a compact Riemann surface of genus $g\geq 2$ and 
$G\subseteq Aut(X)$. If $|G'|$ is odd, then $|G|\leq 30(g-1)$. 
There are infinitely many $g$, necessarily with 
$g\equiv 6\ mod\ 10$, for which this bound is sharp.
\rm
\\ \\
{\bf Proof.} \rm 
In $G/G'$, which is abelian, the elements of odd order form
a unique normal subgroup. So there is a normal subgroup $N$
of $G$ of odd order that contains all elements of odd order.
\hfill$\Box$
\\ \\
{\bf Corollary 7.6.} \it
Let $X$ be a compact Riemann surface of genus $g\geq 2$. 
If $G\subseteq Aut(X)$ with $|G|\equiv 2\ mod\ 4$, then  
$|G|\leq 30(g-1)$. There are infinitely many $g$, necessarily
with $g\equiv 6\ mod\ 10$, for which this bound is sharp.
\rm
\\

\subsection*{8. CLT groups}

\noindent
Another type of groups that has not been treated yet in the
literature in their role as subgroups of $Aut(X)$ are CLT 
groups. The acronym stands for Converse of Lagrange's Theorem.
\par
The precise definition is: A finite group $G$ is CLT if for
every divisor $d$ of $|G|$ there exists a subgroup $H$ of $G$ 
with $|H|=d$.
\par
Every CLT group is solvable, a nontrivial fact. Actually, much
weaker conditions, for example the existence of complements of
the Sylow subgroups, already guarantee solvability 
([H, Theorem 9.3.3] or [R, Theorem 9.1.8]). 
\par
Every supersolvable group is CLT.
Actually, more precisely, a finite group $G$ is supersolvable 
if and only if all its subgroups (including $G$ itself) are CLT 
[Wei, p.13 Theorem 4.1]. 
\par
However, not every solvable group is CLT (smallest counterexample:
$A_4$). And not every CLT group is supersolvable (example $S_4$).
\par
But what makes CLT groups especially annoying to handle: 
A subgroup of a CLT group is not necessarily CLT (example 
$A_4$ in $S_4$). A quotient of a CLT group is not necessarily 
CLT (example $A_4\times C_2$).
\par
To put things into the context of this paper, we now show that 
the last four negative statements persist if one restricts to
groups whose orders are not divisible by primes that are smaller 
than a given prime $p$.
\\ \\
{\bf Example 8.1.} 
Fix an odd prime $p$. By Dirichlet's theorem on primes in 
arithmetic progressions there are infinitely many primes $q$ 
such that $p$ divides $q+1$. Then the finite field $\FF_{q^2}$ 
contains the $p$-th roots of unity, but $\FF_q$ doesn't. 
So over $\FF_q$ the $p$-th cyclotomic polynomial 
$\frac{x^p -1}{x-1}$ splits into irreducible quadratic factors. 
Let a generator of $C_p$ act on the $2$-dimensional vector space 
$\FF_q\oplus\FF_q$ by a $2\times 2$ matrix whose characteristic 
polynomial is one of these irreducible factors. Then there is 
no $1$-dimensional $C_p$-invariant $\FF_q$-subspace. Thus we 
obtain a group
$$H=(C_q\times C_q)\rtimes C_p$$
without subgroups of order $pq$.
So $H$ is not CLT, but obviously it is solvable.
Moreover, $H\times C_q$ is CLT, but not supersolvable.
\rm
\\ \\
We first treat CLT groups of odd order, because for them we
can prove a definitive result.
\\ \\
{\bf Theorem 8.2.} \it
Let $X$ be a compact Riemann surface of genus $g\geq 2$. 
If $G\subseteq Aut(X)$ is a CLT group of odd order, then 
$|G|\leq\frac{21}{2}(g-1)$. The CLT groups of odd order
that reach this bound are exactly the supersolvable groups
described in Theorem 5.1 and Proposition 7.3.
\rm
\\ \\
{\bf Proof.} \rm 
To establish the bound and its sharpness, in view of Theorems 
2.4 and 5.1 it suffices to show that $|G|$ cannot be $15(g-1)$. 
\par
So let's assume $|G|=15(g-1)$, i.e., that $G$ is a quotient of
$\Gamma(0;3,3,5)$. Being a CLT group, $G$ then has a subgroup 
$H$ of index $5$. Left multiplication of $G$ on the cosets 
$G/H$ gives a group homomorphism onto a transitive subgroup 
$U$ of $S_5$. Since $U$ has odd order, the only possibility 
is $U\cong C_5$, contradicting $G/G'\cong C_3$.
\par
Now let $G$ be a CLT group of odd order $|G|=\frac{21}{2}(g-1)$.
By Proposition 7.3 we still have to show $|G|=3\cdot 7^n$.
By exactly the same proof as above we see that $5$ cannot
divide $|G|$. 
\par
Now assume that $9$ divides $|G|$. Then $G$ has a subgroup
of index $9$. As above, we get a homomorphism from $G$ onto 
a transitive subgroup $U$ of $S_9$ of order $3^{\mu}7^{\nu}$ 
with $2\leq\mu\leq 4$ and $\nu\leq 1$. In any case, the Sylow 
$7$-subgroup of $U$ is either trivial or normal, leading to 
a quotient of $G$ of order $3^{\mu}$, and ultimately to 
a quotient of $\Gamma(0;3,3,7)$ of order $9$. This shows
that $9$ cannot divide $|G|$.
\par
Finally, let $p$ be the smallest prime divisor of $|G|$ that 
is bigger than $7$. The following proof that $p$ cannot exist 
uses the same trick that was used in the proof of 
[Sch1, Proposition 5.3], following a very helpful suggestion 
by the referee of that paper. 
\par
Since $G$ is a CLT group, it has a subgroup of index $p$. 
Let $N$ be its core. Then $[G:N]$ can only be divisible
by $p$ (once), by $3$ (once), and by a power of $7$. Let $U$
be the quotient of $G/N$ by its biggest normal $7$-subgroup.
Since the only abelian quotient of $\Gamma(0;3,3,7)$ is $C_3$,
we have $[U:U']=3$ and $U'$ is a quotient of $\Gamma(0;7,7,7)$.
So we obtain a sequence of normal subgroups 
$$U\triangleright U' \triangleright P\triangleright I$$
where $P$ ($\cong C_p$) is the Sylow $p$-subgroup (which is 
normal in $U$) and $U'/P$ is isomorphic to the not necessarily 
abelian Sylow $7$-subgroup $S$ of $U$. Thus $U'\cong P\rtimes S$.
Now let $C$ be the centralizer of $P$ in $U$. As $U/C$ can be 
embedded in the automorphism group of $P$, which is cyclic, we
see that $C$ contains $U'$, and hence that $P$ is central in 
$U'$. So $U'$ is a direct product of $P$ and $S$. In particular, 
$P$ is a quotient of $U'$ and hence of $\Gamma(0;7,7,7)$, giving 
the desired contradiction.
\hfill$\Box$
\\ \\
As CLT groups are solvable, in general we have the bound 
$|G|\leq 48(g-1)$. We refine this a little bit. For the 
proof we recall that a Hall subgroup $H$ of a finite group 
$G$ is a subgroup whose order $|H|$ is relatively prime 
to its index $[G:H]$. 
\\ \\
{\bf Lemma 8.3.} \it
A CLT group $G$ that reaches the bound $|G|=48(g-1)$ 
must necessarily be a $(2,3)$-group. 
\rm
\\ \\
\noindent
{\bf Proof.} \rm 
If $G\subseteq Aut(X)$ is a CLT group of order $48(g-1)$,
let $p$ be the smallest prime divisor of $|G|$ that is bigger 
than $3$. The following proof that $p$ cannot exist uses the 
same trick as in the proof of Theorem 8.2, due to the referee
of the paper [Sch1].
\par
Since $G$ is a CLT group, it has a subgroup of index $p$. 
Let $N$ be its core. Then $[G:N]$ can only be divisible
by $p$ (once) and by powers of $2$ and $3$. Let $U$
be the quotient of $G/N$ by its biggest normal subgroup
of order prime to $p$.
Since the only abelian quotient of $\Gamma(0;2,3,8)$ is $C_2$,
we have $[U:U']=2$ and $U'$ is a quotient of $\Gamma(0;3,3,4)$.
So we obtain a sequence of normal subgroups 
$$U\triangleright U' \triangleright P\triangleright I$$
where $P$ ($\cong C_p$) is the Sylow $p$-subgroup (which is 
normal in $U$) and $U'/P$ is isomorphic to a not necessarily 
abelian Hall subgroup $H$ of $U'$ (see [H, Theorem 9.3.1]). 
Thus $U'\cong P\rtimes H$.
Now let $C$ be the centralizer of $P$ in $U$. As $U/C$ can be 
embedded in the automorphism group of $P$, which is cyclic, we
see that $C$ contains $U'$, and hence that $P$ is central in 
$U'$. So $U'$ is a direct product of $P$ and $H$. In particular, 
$P$ is a quotient of $U'$ and hence of $\Gamma(0;3,3,4)$, giving 
the desired contradiction.
\hfill$\Box$
\\ \\
{\bf Example 8.4.} 
The group of order $48$ acting on a Riemann surface 
$X$ of genus $2$ is 
$$G\cong GL_2(\FF_3)\cong Q_8 \rtimes S_3$$ 
where $Q_8$ is the quaternion group of order $8$. Let 
$Z$($\cong C_2$) be the center of $Q_8$ and $C_3$ the subgroup 
of $S_3$. Then it is easy to see that $G$ is a CLT group. 
Subgroups of order $6$, $12$ and $24$ are $Z\times  C_3$, 
$Z\rtimes S_3$ and $Q_8 \rtimes C_3$. For the other orders 
the existence of a subgroup is clear from the Sylow Theorems.
\par
By similar arguments one can show that the group of order 
$96$ that acts on a Riemann surface of genus $3$, namely
$$G\cong (C_4 \times C_4)\rtimes S_3$$
also is a CLT group.
\\ \\
{\bf Remark 8.5.}
Originally, our intention was to prove the existence of 
infinitely many CLT groups of order $48(g-1)$ by applying
Macbeath's approach to one of the groups from Example 8.4. 
\par
Let for example $G$ be the group of order $48$ acting on
a Riemann surface $X$ of genus $2$. For every $\mu\geq 1$
there exists a totally unramified Galois covering $Y$ of $X$ 
with $Gal(Y/X)\cong(C_{2^{\mu}})^4$. Then $Aut(Y)$ contains 
a subgroup $\widetilde{G}$ with
$Gal(Y/X)\triangleleft\widetilde{G}$ and 
$\widetilde{G}/Gal(Y/X)\cong G$. 
\par
{\bf If} one manages to show that
$$\widetilde{G}\cong Gal(Y/X)\rtimes G,$$
then $\widetilde{G}$ is easily seen to be a CLT group
by the following argument. Let $d$ be a divisor of 
$|\widetilde{G}|=3\cdot 2^{4\mu +4}$.
Choose the biggest $\nu\leq\mu$ such that $2^{4\nu}$ divides
$d$. Then $d=2^{4\nu}d'$ with $d'|48$. We pick the subgroup
$(C_{2^{\nu}})^4$ of $Gal(Y/X)$, which is characteristic in
$Gal(Y/X)$ and hence stable under the action of $G$. And 
we also pick a subgroup $U$ of order $d'$ of the CLT group 
$G$. Then $(C_{2^{\nu}})^4 \rtimes U$ is the desired subgroup 
of $\widetilde{G}$ of order $d$.
\par
However, we have not been able to show the semi-direct 
product structure or something comparable that would make
our attempt work.
\par
On the other hand, it might equally well be that the two 
groups in Example 8.4 are indeed the only CLT groups of 
order $48(g-1)$.
\\ \\
As Theorem 8.2 and Lemma 8.3 are probably the first results 
in the literature that deal with CLT groups as automorphism 
groups of Riemann surfaces, we elaborate a bit more on this 
encounter. For that we need the following tools. 
\\ \\
{\bf Theorem 8.6.} \it
Let $G$ be a CLT group. Then
\begin{itemize}
\item[(a)] Every Hall subgroup $H$ of $G$ is also CLT.
\item[(b)] $[G:G']$ is divisible by the smallest prime divisor of $|G|$.
\item[(c)] If $4$ divides $|G|$, then $G$ has a quotient $C_4$ or 
$C_2 \times C_2$ or $C_6$ or $S_4$.
\end{itemize}

\rm
\noindent
{\bf Proof.} \rm 
(a) Let $d$ be a divisor of $m=|H|$. Being a CLT group, $G$ has 
a subgroup $U$ of order $d$. But by [H, Theorem 9.3.1] any subgroup
of order $d$ is contained in a Hall subgroup of order $m$, and all
Hall subgroups of order $m$ are conjugate. So $H$ contains a conjugate 
of $U$.
\par
(b) Let $p$ be the smallest prime dividing $|G|$. Being a CLT group,
$G$ has a subgroup $U$ of index $p$, and since $p$ is the smallest 
prime dividing $|G|$, this subgroup is normal; so $G'\subseteq U$.
\par
(c) Being a CLT group, $G$ has a subgroup $U$ of index $4$. Let $N$
be its core. Then $G/N$ is isomorphic to one of $C_4$, $C_2 \times C_2$,
$D_4$, $A_4$, $S_4$. The first three imply the existence of a quotient
of order $4$. And $A_4$ has a quotient of order $3$, which together 
with part (b) shows that $G/G'$ has a quotient $C_6$.
\hfill$\Box$
\\ \\
{\bf Lemma 8.7.} \it
There are no CLT groups among the groups of orders $40(g-1)$, $36(g-1)$, 
$30(g-1)$, $21(g-1)$, $20(g-1)$, $\frac{56}{3}(g-1)$.
\rm
\\ \\
{\bf Proof.} \rm 
Compare Table 2. 
\par
The case $|G|=21(g-1)$ is the most tricky one. Then $G$ is a quotient 
of $\Gamma(0;2,3,14)$. Let $G$ also be a CLT group. Then $4$ cannot
divide $|G|$, because $\Gamma(0;2,3,14)$ has none of the quotients 
listed in Theorem 8.6 (c). On the other hand, $C_2$ is the only
abelian quotient of $\Gamma(0;2,3,14)$. So $G'$ is a Hall subgroup
of $G$, and hence must also be CLT by Theorem 8.6 (a). Consequently,
$|G'|=3\cdot 7^n$ by Theorem 8.2. 
Being a CLT group, $G$ has a subgroup of index $7$. Its core is 
a normal subgroup $N$ of $G$ of index $42$. Here we are using again
(twice) that if $G$ is solvable and $[G:G']$ a prime, then $G'$
contains all normal subgroups of $G$. Denote $G/N$ by $H$. Since the 
only abelian quotients of $\Gamma(0;2,3,14)$ and of $\Gamma(0;3,3,7)$
are $C_2$ resp. $C_3$, we get $H/H'\cong C_2$, $H'/H''\cong C_3$ and
$H''/H'''\cong C_7$. But this contradicts [H, Theorem 9.4.2], which
says that $H'/H''$ and $H''/H'''$ cannot both be cyclic.
\par
The case $\Gamma(0;2,3,10)$, i.e., $|G|=30(g-1)$ is discarded 
similarly. Theorem 8.6 (c) shows that $4$ cannot divide $|G|$; 
then Theorem 8.6 (a) implies the existence of an odd order CLT 
group $H$ with $|H|=15(g-1)$, in contradiction to Theorem 8.2.
\par
The other cases are also quickly excluded by Theorem 8.6. 
The only quotient of $\Gamma(0;2,4,5)$ of order prime 
to $5$ is $C_2$; and $C_2$ is also the only quotient 
of $\Gamma(0;2,4,7)$ of order prime to $7$. 
The groups $\Gamma(0;2,3,9)$, $\Gamma(0;2,3,15)$ and 
$\Gamma(0;2,5,5)$ do not even have a quotient $C_2$. 
\hfill$\Box$
\\ \\
For the same reasons the triangle group $\Gamma(0;3,3,4)$ can 
never furnish a CLT group. However, although $\Gamma(0;3,3,4)$
is the commutator group of $\Gamma(0;2,3,8)$, this does not 
exclude the possibility of CLT groups of order $48(g-1)$.
The reason simply is that the commutator group of a CLT group, 
even if it has index $2$, is not necessarily a CLT group, as is
already shown by the standard example $S_4$. 
\\ \\
{\bf Corollary 8.8.} \it
The sharp bound for CLT groups $G$ whose order is not divisible 
by $8$ is the same as for supersolvable groups, namely 
$|G|\leq 18(g-1)$ if $g\geq 3$.
\rm
\\ \\
{\bf Proof.} \rm 
From Table 2 we see that the orders bigger than $18(g-1)$ that 
were not already excluded by Lemma 8.7 are all divisible by $8$.
\hfill$\Box$
\\

\subsection*{\hspace*{10.5em} References}
\begin{itemize}

\item[{[BMT]}] B.~Borror, A.~Morris, M.~Tarr: 
\rm The strong symmetric genus spectrum of abelian groups, 
\it Rose-Hulman Undergrad. Math. J. \bf 15 \rm (2014), 115-130 

\item[{[B]}] T.~Breuer: \it Characters and Automorphism Groups 
of Compact Riemann Surfaces, \rm LMS Lecture Notes 280, Cambridge 
University Press, Cambridge, 2000

\item[{[BC]}] E.~Bujalance and M.~Conder:
\rm On cyclic groups of automorphisms of Riemann surfaces,
\it J. London Math. Soc. \bf 59 \rm (1999), 573-584 

\item[{[Ch]}] B.~P.~Chetiya: On genuses of compact Riemann 
surfaces admitting solvable automorphism groups, 
\it Indian J. Pure Appl. Math. \bf 12 \rm (1981), 1312-1318 

\item[{[ChP]}] B.~P.~Chetiya and K.~Patra: On metabelian 
groups of automorphisms of compact Riemann surfaces, 
\it J. London Math. Soc. \bf 33 \rm (1986), 467-472

\item[{[FT]}] W.~Feit and J.~G.~Thompson: Solvability of groups 
of odd order, \it Pacific J. Math. \bf 13 \rm (1963), 775-1029 

\item[{[G1]}] G.~Gromadzki: Maximal groups of automorphisms of
compact Riemann surfaces in various classes of finite groups,
\it Rev. Real Acad. Cienc. Exact. Fis. Natur. Madrid \bf 82 no. 2
\rm (1988), 267-275

\item[{[G2]}] G.~Gromadzki: On soluble groups of automorphism of 
Riemann surfaces, \it Canad. Math. Bull. \bf 34 \rm (1991), 67-73

\item[{[G3]}] G.~Gromadzki: Metabelian groups acting on 
compact Riemann surfaces, \it Rev. Mat. Univ. Complut. Madrid 
\bf 8 no. 2 \rm (1995), 293-305

\item[{[GMl]}] G.~Gromadzki and C.~Maclachlan: Supersoluble
groups of automorphisms of compact Riemann surfaces, 
\it Glasgow Math. J. \bf 31 \rm (1989), 321-327

\item[{[KS]}] S.~Kallel and D.~Sjerve: \rm On the group of 
automorphisms of cyclic covers of the Riemann sphere,
\it Math. Proc. Camb. Phil. Soc. \bf 138 \rm (2005), 267-287 

\item[{[H]}] M.~Hall: \it The Theory of Groups, \rm Macmillan,
New York, 1959 

\item[{[Ha]}] W.~J.~Harvey: \rm Cyclic groups of automorphisms 
of a compact Riemann surface, \it Quart. J. Math. Oxford (2), 
\bf 17 \rm (1966), 86-97

\item[{[Mb]}] A.~M.~Macbeath: On a theorem of Hurwitz,
\it Proc. Glasgow Math. Assoc. \bf 5 \rm (1961), 90-96 

\item[{[Ml]}] C.~Maclachlan: \rm Abelian groups of automorphisms 
of compact Riemann surfaces, \it Proc. London Math. Soc. \bf 15
\rm (1965), 699-712

\item[{[MZ1]}] C.~May and J.~Zimmerman: \rm The symmetric genus 
of metacyclic groups, \it Topol. Appl. \bf 66 \rm (1995), 101-115

\item[{[MZ2]}] C.~May and J.~Zimmerman: The symmetric genus
of groups of odd order, \it Houston J. Math. \bf 34 no. 2 
\rm (2008), 319-338

\item[{[N]}] K.~Nakagawa: On the orders of automorphisms of 
a closed Riemann surface, \it Pacific J. Math. \bf 115 no. 2 
\rm (1984), 435-443 

\item[{[PR]}] R.~Pignatelli and C.~Raso: Riemann surfaces
with a quasilarge abelian group of automorphisms, \it Matematiche 
\bf 66 no. 2 \rm (2011), 77-90 

\item[{[R]}] D.~J.~S.~Robinson: \it A Course in the Theory of Groups,
\rm Springer GTM 80, New York - Berlin, 1982

\item[{[Sch1]}] A.~Schweizer: On the exponent of the automorphism 
group of a compact Riemann surface, \it Arch. Math. (Basel) 
\bf 107 \rm (2016), 329-340 

\item[{[Sch2]}] A.~Schweizer: Metacyclic groups as automorphism
groups of compact Riemann surfaces, \it Geom. Dedicata \rm
(2017), DOI: 10.1007/s10711-017-0239-8 

\item[{[W1]}] A.~Weaver: Genus spectra for split metacyclic 
groups, \it Glasgow Math. J. \bf 43 \rm (2001), 209-218 

\item[{[W2]}] A.~Weaver: Odd-order groups actions on surfaces,
\it J. Ramanujan Math. Soc. \bf 18\rm , no. 3 (2003), 211-220

\item[{[Wei]}] M.~Weinstein (editor): \it Between nilpotent and 
solvable, \rm Polygonal Publ. House, Washington, N. J., 1982 

\item[{[Z1]}] R.~Zomorrodian: Nilpotent automorphism groups of 
Riemann surfaces, \it Trans. Amer. Math. Soc. \bf 288 no. 1 
\rm (1985), 241-255

\item[{[Z2]}] R.~Zomorrodian: Classification of $p$-groups of 
automorphisms of Riemann surfaces and their lower central series,
\it Glasgow Math. J. \bf 29 \rm (1987), 237-244

\item[{[Z3]}] R.~Zomorrodian: Bounds for the order of supersoluble
automorphism groups of Riemann surfaces, \it Proc. Amer. Math. Soc.
\bf 108 no. 3 \rm (1990), 587-600

\item[{[Z4]}] R.~Zomorrodian: On a theorem of supersoluble 
automorphism groups, \it Proc. Amer. Math. Soc.
\bf 131 no. 9 \rm (2003), 2711-2713

\end{itemize}

\end{document}